\documentclass[titlepage,twoside,12pt]{article}
\usepackage{amssymb}
\usepackage{amsfonts}
\begin{document}

{\bf \Large Heisenberg Uncertainty in Reduced Power \\ \\ Algebras} \\ \\

{\bf Elem\'{e}r E Rosinger} \\ \\
Department of Mathematics \\
and Applied Mathematics \\
University of Pretoria \\
Pretoria \\
0002 South Africa \\
eerosinger@hotmail.com \\ \\

{\bf Abstract} \\

The Heisenberg uncertainty relation is known to be obtainable by a purely mathematical argument. Based on that fact, here it is shown that the Heisenberg
uncertainty relation remains valid when Quantum Mechanics is re-formulated within far wider frameworks of {\it scalars}, namely, within one or the other
of the infinitely many {\it reduced power algebras} which can replace the usual real numbers $\mathbb{R}$, or complex numbers $\mathbb{C}$. A major
advantage of such a re-formulation is, among others, the disappearance of the well known and hard to deal with problem of the so called "infinities in
Physics". The use of reduced power algebras also opens up a foundational question about the role, and in fact, about the very meaning and existence, of
fundamental constants in Physics, such as Planck's constant $h$. A role, meaning, and existence which may, or on the contrary, may not be so objective as
to be independent of the scalars used, be they the usual real numbers $\mathbb{R}$, complex numbers $\mathbb{C}$, or scalars given by any of the
infinitely many reduced power algebras, algebras which can so easily be constructed and used. \\ \\

{\bf 1. Preliminaries} \\

A remarkable feature of the Heisenberg uncertainty relation is that it can be obtained following a purely mathematical argument of a rather simple statistical nature, [1, pp. 67-70]. Based on that fact, here we shall show that the Heisenberg uncertainty relation remains valid when Quantum Mechanics is re-formulated within what appears to be a
far more wide and appropriate framework of {\it scalars}, namely, those in any of the infinitely many algebras which belong to the class of {\it reduced power algebras}, [2-9]. \\

As argued in [2-9], there are a number of important advantages in re-formulating the whole of Physics in terms of scalars given by reduced power algebras. Related to Quantum Theory, and specifically, to Quantum Field Theory, one of the major advantages of such a re-formulation is the complete and automatic disappearance of the well known and hard to deal with problem of the so called "infinities in Physics". \\

In this regard, let us recall that in [2-9] it was shown how to construct in a simple way as {\it reduced powers} a large class of algebras which extend the field $\mathbb{R}$ of usual real numbers, or alternatively, the field $\mathbb{C}$ of usual complex numbers. A remarkable feature of these {\it reduced power algebras} is that they contain {\it infinitesimal}, as well as {\it infinitely large} elements, consequently, these algebras are {\it non-Archimedean}. Some of these algebras are in fact fields. Also, among them is the field $^*\mathbb{R}$ of nonstandard real numbers. \\

In [2-9] it was suggested and argued that much of present day Physics should be re-formulated in terms of such reduced power algebras, one of the main reasons for that being the considerably increased richness and complexity of their non-Archimedean self-similar mathematical structures, as opposed to the much simpler structures imposed on Physics by the Archimedean field $\mathbb{R}$ of the usual real numbers, and by the structures built upon it, such as the field $\mathbb{C}$ of complex numbers, various finite or infinite dimensional manifolds, Hilbert spaces, and so on. And in this regard, it was argued that one of the main advantages of such a re-formulation would be the automatic disappearance of the difficulties related to the so called "infinities in Physics", as a result of the presence of {\it infinitely large} elements in the reduced power algebras. \\

It was also pointed out that the present limitation to the exclusive use of scalars, vectors, etc., which belong to Archimedean mathematical structures is the result not of absolutely any kind of conscious and competent choice in Physics, but on the contrary, of the millennia long perpetuation of a mere historical accident, namely that in ancient Egypt the development of Geometry chose the Archimedean route, due to specific practical needs at those times, needs hardly at all related to those of modern Physics. \\

However, as not seldom happens in human affairs, accidentally acquired habits can become second nature. This may explain, even if not excuse as well,
why for more than two millennia by now we have been so happily wallowing in the ancient Egyptian bondage, or rather slavery of Archimedean space-time structures ... \\

For convenience, we shall recall in a particular case the construction, [2-9], of {\it reduced power algebras}. Given any {\it filter} ${\cal F}$ on $\mathbb{N}$, we define \\

(1.1)~~~ $ {\cal I}_{\cal F} = \{~ v = ( v_n )_{n \in \mathbb{N}} \in \mathbb{R}^\mathbb{N} ~~|~~
                                       \{~ n \in \mathbb{N} ~|~ v_n = 0 ~\} \in {\cal F} ~\} $ \\

which is a {\it proper ideal} in the algebra $\mathbb{R}^\mathbb{N}$. Thus we obtain the {\it reduced power algebra} associated to ${\cal F}$ as the quotient algebra \\

(1.2)~~~ $ \mathbb{R}_{\cal F} = \mathbb{R}^\mathbb{N} / {\cal I}_{\cal F} $ \\

Furthermore, this algebra which is commutative, is also a strict extension of the field $\mathbb{R}$ of the usual real numbers, according to the embedding of algebras \\

(1.3)~~~ $ \mathbb{R} \ni x \longmapsto ( x, x, x, \ldots ) + {\cal I}_{\cal F} \in \mathbb{R}_{\cal F} = \mathbb{R}^\mathbb{N} / {\cal I}_{\cal F} $ \\

In a similar manner one can obtain reduced power algebras extending the field $\mathbb{C}$ of the usual complex numbers. Namely, let us denote by \\

(1.4)~~~ $ {\cal J}_{\cal F} = \{~ w = ( w_n )_{n \in \mathbb{N}} \in \mathbb{C}^\mathbb{N} ~~|~~
                                       \{~ n \in \mathbb{N} ~|~ w_n = 0 ~\} \in {\cal F} ~\} $ \\

which is a {\it proper ideal} in the algebra $\mathbb{C}^\mathbb{N}$. Thus we obtain the {\it reduced power algebra} associated to ${\cal F}$ as the quotient algebra \\

(1.5)~~~ $ \mathbb{C}_{\cal F} = \mathbb{C}^\mathbb{N} / {\cal J}_{\cal F} $ \\

Furthermore, this algebra which is commutative, is also a strict extension of the field $\mathbb{C}$ of the usual complex numbers, according to the embedding of algebras \\

(1.6)~~~ $ \mathbb{C} \ni z \longmapsto ( z, z, z, \ldots ) + {\cal I}_{\cal F} \in \mathbb{C}_{\cal F} =
                             \mathbb{C}^\mathbb{N} / {\cal J}_{\cal F} $ \\

We now establish a natural connection between the algebras $\mathbb{R}_{\cal F}$ and $\mathbb{C}_{\cal F}$. \\

In this regard, we note the following connection between the ideals ${\cal I}_{\cal F}$ and ${\cal J}_{\cal F}$. Namely \\

(1.7)~~~ $ \begin{array}{l}
                  w = ( w_n = u_n + i v_n )_{n \in \mathbb{N}} \in {\cal J}_{\cal F} ~~\Longleftrightarrow \\ \\
                  ~~~~~~~\Longleftrightarrow~~ u = ( u_n )_{n \in \mathbb{N}},~
                             v = ( v_n )_{n \in \mathbb{N}} \in {\cal I}_{\cal F}
           \end{array} $ \\

where $u_n, v_n \in \mathbb{R}$. It follows that we have the algebra homomorphisms \\

(1.8)~~~ $ \begin{array}{l}
              Re : \mathbb{C}_{\cal F} \ni  w = ( w_n = u_n + i v_n )_{n \in \mathbb{N}} + {\cal J}_{\cal F} \longmapsto \\ \\
              ~~~~~~ \longmapsto
                 u = ( u_n )_{n \in \mathbb{N}} + {\cal I}_{\cal F} \in \mathbb{R}_{\cal F}
           \end{array} $ \\ \\

(1.9)~~~ $ \begin{array}{l}
               Im : \mathbb{C}_{\cal F} \ni  w = ( w_n = u_n + i v_n )_{n \in \mathbb{N}} + {\cal J}_{\cal F} \longmapsto \\ \\
               ~~~~~~ \longmapsto
                 v = ( v_n )_{n \in \mathbb{N}} + {\cal I}_{\cal F} \in \mathbb{R}_{\cal F}
           \end{array} $ \\

as well as the algebra embeddings \\

(1.10)~~~ $ \mathbb{R}_{\cal F} \ni u = ( u_n )_{n \in \mathbb{N}} + {\cal I}_{\cal F} \longmapsto
                            u = ( u_n )_{n \in \mathbb{N}} + {\cal J}_{\cal F} \in \mathbb{C}_{\cal F} $ \\

(1.11)~~~ $ \mathbb{R}_{\cal F} \ni v = ( v_n )_{n \in \mathbb{N}} + {\cal I}_{\cal F} \longmapsto
                            i v = ( i v_n )_{n \in \mathbb{N}} + {\cal J}_{\cal F} \in \mathbb{C}_{\cal F} $ \\

Let us also define the surjective linear mapping \\

(1.12)~~~ $ \begin{array}{l}
                 \mathbb{C}_{\cal F} \ni  w = ( w_n = u_n + i v_n )_{n \in \mathbb{N}} + {\cal J}_{\cal F} \longmapsto \\ \\
                 ~~~~~~ \longmapsto
                 \overline{w} = ( \overline{w_n} = u_n - i v_n )_{n \in \mathbb{N}} + {\cal J}_{\cal F} \in \mathbb{C}_{\cal F}
             \end{array} $ \\

As a consequence, we obtain \\

(1.13)~~~ $ w = ( w_n = u_n + i v_n )_{n \in \mathbb{N}} + {\cal J}_{\cal F} \in \mathbb{C}_{\cal F},~~~
                  \overline{w} = w
           ~~~\Longrightarrow~~~ w \in \mathbb{R}_{\cal F} $ \\

Lastly, we can define the {\it absolute value} on $\mathbb{C}_{\cal F}$, by the mapping \\

(1.14)~~~ $ \begin{array}{l}
                \mathbb{C}_{\cal F} \ni z = ( w_n = u_n + i v_n )_{n \in \mathbb{N}} + {\cal J}_{\cal F}
                                 \longmapsto \\ \\
                ~~~~~~ \longmapsto
                   | z | = ( | w_n | = \sqrt ( u^2_n + v^2_n ) )_{n \in \mathbb{N}} + {\cal I}_{\cal F} \in \mathbb{R}_{\cal F}
            \end{array} $ \\

Let us denote \\

(1.15)~~~ $ \mathbb{R}^+_{\cal F} = \{~ u = ( u_n )_{n \in \mathbb{N}} + {\cal I}_{\cal F} \in \mathbb{R}_{\cal F}
                     ~~|~~ \{~ n \in \mathbb{N} ~|~ u_n \geq 0 ~\} \in {\cal F} ~\} $ \\

then we obtain the surjective mapping \\

(1.16)~~~ $  \mathbb{C}_{\cal F} \ni z \longmapsto | z | \in \mathbb{R}^+_{\cal F} $ \\

and for $z \in \mathbb{C}_{\cal F}$, we have \\

(1.17)~~~ $ | z | = 0 ~~\Longleftrightarrow~~ z = 0 $ \\

Now, in view of (1.8), (1.9), (1.14), we have for $z \in \mathbb{C}_{\cal F}$ the relations \\

(1.18)~~~ $ |\, Re \, z \,|,~~ |\, Im \, z \,| ~\leq~ |\, z \,| $ \\

For further convenience, we shall consider a {\it quantum configuration space} which is one dimensional. \\

Here however, there are two ways to proceed. \\

The simpler one is to model the one dimensional configuration space by the usual $\mathbb{R}$, in which case the {\it wave functions} will be given by \\

(1.19)~~~ $ \psi : \mathbb{R} \longrightarrow \mathbb{C}_{\cal F} $ \\

This means that the only difference with the usual quantum mechanical setup is that, this time, the wave functions can take values in the reduced power algebra extension $\mathbb{C}_{\cal F}$ of $\mathbb{C}$. \\

Alternatively, one can be more consistent in the re-formulation of Physics in terms of reduced power algebras, and
model not only the values of the wave functions, but also their one dimensional configuration space variables with the reduced power algebra $\mathbb{R}_{\cal F}$ which is an extension of $\mathbb{R}$. Thus in this second case, the wave functions would be \\

(1.20)~~~ $ \psi : \mathbb{R}_{\,\cal U} \longrightarrow \mathbb{C}_{\cal F} $ \\

where ${\cal U}$ is an ultrafilter on $\mathbb{N}$. \\

The first of these two alternatives will be developed in the next section. The second and yet more general alternative is treated elsewhere. \\ \\

{\bf 2. An Extension of the Heisenberg Uncertainty} \\

In order to avoid unnecessary technical complications concerning the integrations related to wave functions $\psi$ in (1.19), we shall only consider
those of them which are of the following particular {\it step function} type \\

(2.1)~~~ $ \psi = \sum_{1 \leq h \leq m} \gamma_h H_h $ \\

where $m \geq 1$, and $\gamma_h \in \mathbb{C}_{\cal F}$, while $H_h : \mathbb{R} \longrightarrow \{ 0, 1 \}$ are step functions such that $H_h ( x ) =
1$, when $x \in I_h$, and $H_h ( x ) = 0$, when $x \notin I_h$. Here $I_h = [ a_{h-1}, a_h) \subset \mathbb{R}$ are usual intervals, where $\infty \leq
a_0 < a_1 < a_2 \leq \ldots < a_m \leq \infty$ are usual real numbers, with $a_0$ and $a_m$ possibly minus and plus infinity, respectively. \\

The set of all such functions wave functions $\psi$ in (2.1) is denoted by \\

(2.2)~~~ $ {\cal S}_{\cal F} ( \mathbb{R} ) $ \\

and this set replaces the usual Hilbert space $L^2 ( \mathbb{R} )$ of complex valued square integrable wave functions $\psi$ defined on the
configuration space $\mathbb{R}$ of a one dimensional quantum system. \\

We note here one of the major {\it advantages} in the above use of reduced power algebras. Namely, each of the wave functions $\psi \in {\cal S}_{\cal
F} ( \mathbb{R} )$ can easily be integrated on the whole of $\mathbb{R}$ regardless of the possible infinite length of some of the intervals $I_h$, or
of the infinite value of some of the $\gamma_h$. Indeed, owing to the reduced power algebra structure of $\mathbb{C}_{\cal F}$, one simply obtains in
$\mathbb{C}_{\cal F}$ the algebraically perfectly well defined value \\

(2.3)~~~ $ \int_{\mathbb{R}} \psi ( x ) dx =
                 \sum_{1 \leq h \leq m} ( a_h - a_{h-1} ) \gamma_h \in \mathbb{C}_{\cal F} $ \\

which is always a well defined element in $\mathbb{C}_{\cal F}$, and thus available in a correct and rigorous manner for all the algebraic operations in
the algebra $\mathbb{C}_{\cal F}$, even if that value may turn out to be an infinitesimal, finite, or an infinitely large element in $\mathbb{C}_{\cal
F}$. In this way there is no need to impose any integrability type conditions on the wave functions $\psi \in {\cal S}_{\cal F} ( \mathbb{R} )$, much
unlike in the usual case, where the square integrability condition $\int_{\mathbb{R}} | \psi ( x ) |^2 dx < \infty$ is required, since within the usual
Archimedean framework of $\mathbb{C}$, or $\mathbb{R}$, one cannot perform most of the usual algebraic operations with infinitely large quantities. \\

Now it follows easily that ${\cal S}_{\cal F} ( \mathbb{R} )$ is a vector space over $\mathbb{C}$, and in fact, it is a commutative algebra over $\mathbb{C}$. Furthermore, one can define on it the extension of the usual scalar product given by \\

(2.4)~~~ $ < \psi, \chi > ~=~ \int_{\mathbb{R}} \overline{\psi ( x )} \chi ( x )dx \in \mathbb{C}_{\cal F} $ \\

for all $\psi, \chi \in {\cal S}_{\cal F} ( \mathbb{R} )$. This extended scalar product has the following properties : \\

(2.5)~~~ It is linear over $\mathbb{C}_{\cal F}$, therefore also over $\mathbb{C}$, in the second \\
         \hspace*{1.3cm} argument. \\

(2.6)~~~ $  < \chi, \psi > ~=~ \overline{< \psi, \chi >},~~~ \psi, \chi \in {\cal S}_{\cal F} ( \mathbb{R} ) $ \\

(2.7)~~~ $ < \psi, \psi > \,\, \in \mathbb{R}^+_{\cal F},~~~ \psi \in {\cal S}_{\cal F} ( \mathbb{R} ) $ \\

and for $\psi \in {\cal S}_{\cal F} ( \mathbb{R} )$, one has \\

(2.8)~~~ $ < \psi, \psi > ~=~ 0 ~~\Longleftrightarrow~~ \psi = 0 \in {\cal S}_{\cal F} ( \mathbb{R} ) $ \\

Also, we have the extension of the classical Schwartz inequality \\

(2.9)~~~ $ |\, < \psi, \chi > \,\,| ~\leq~ < \psi, \psi >^{1/2} \,\, < \chi, \chi >^{1/2},~~~
                                  \psi, \chi \in {\cal S}_{\cal F} ( \mathbb{R} ) $ \\

Laastly, we can consider the set \\

(2.10)~~~ $ {\cal L} ( {\cal S}_{\cal F} ( \mathbb{R} ) ) $ \\

of all {\it linear operators} $A : {\cal S}_{\cal F} ( \mathbb{R} ) \longrightarrow {\cal S}_{\cal F} ( \mathbb{R} )$. Such an operator will be called {\it Hermitian}, if and only if \\

(2.11)~~~ $ < A \psi, \chi > ~=~ < \psi, A \chi >,~~~ \psi, \chi \in {\cal S}_{\cal F} ( \mathbb{R} )$ \\

and we denote by \\

(2.12)~~~ $ {\cal H} ( {\cal S}_{\cal F} ( \mathbb{R} ) ) $ \\

the set of all such Hermitian operators. \\

With these preparations, we can now proceed to obtain the Heisenberg uncertainty relation for arbitrary operators $A, B \in {\cal H} ( {\cal S}_{\cal F} ( \mathbb{R} ) )$. \\

Given $A \in {\cal L} ( {\cal S}_{\cal F} ( \mathbb{R} ) )$ and $\psi \in {\cal S}_{\cal F} ( \mathbb{R} )$, we
denote \\

(2.13)~~~ $ < A >_\psi ~=~ < \psi, A \psi > \,\, \in \mathbb{C}_{\cal F} $ \\

and call it the {\it expectation value} of $A$ in the state $\psi$. Further, we denote \\

(2.14)~~~ $ \Delta_{\psi} A = (\, < A^2 >_\psi \,-\, ( < A >_\psi )^2 \,)^{1 / 2} $ \\

and call it the {\it uncertainty} of $A$ in the state $\psi$. \\

{\bf Theorem 2.1. (Extended Heisenberg Uncertainty Relation)} \\

Given $A, B \in {\cal H} ( {\cal S}_{\cal F} ( \mathbb{R} ) )$ and $\psi \in {\cal S}_{\cal F} ( \mathbb{R} )$ such that $ < \psi, \psi > \,= 1$, then we have \\

(2.15)~~~ $ \Delta_\psi A \,\, \Delta_\psi B ~\geq~ | < [ A, B ] >_\psi | / 2 $ \\

where $[ A, B ] = A B - B A$. \\

{\bf Proof.} \\

We start with \\

{\bf Lemma 2.1.} \\

Let $A \in {\cal H} ( {\cal S}_{\cal F} ( \mathbb{R} ) )$, $\psi \in {\cal S}_{\cal F} ( \mathbb{R} )$, then \\

(2.16)~~~ $ < A >_\psi \,\in\, \mathbb{R}_{\cal F} $ \\

{\bf Proof.} \\

We have in view of (2.13), (2.11) \\

$~~~~~~ < A >_\psi ~=~ < \psi, A \psi > ~=~ < \psi, A \psi > $ \\

thus (2.4), (2.3) imply \\

$~~~~~~ < \psi, A \psi > ~=~ \overline{< \psi, A \psi >} $ \\

hence (1.13) completes the proof.

\hfill $\Box$ \\

Let us now denote \\

(2.17)~~~ $ A_1 = A \,\,-\, < A >_\psi I,~~~ B_1 = B \,\,-\, < B >_\psi I $ \\

where $I \in {\cal H} ( {\cal S}_{\cal F} ( \mathbb{R} ) )$ is the identity operator. Then \\

(2.18)~~~ $ A_1, B_1 \in {\cal H} ( {\cal S}_{\cal F} ( \mathbb{R} ) ) $ \\

Indeed, let $\eta, \chi \in {\cal S}_{\cal F} ( \mathbb{R} )$, then in view of (2.11) and Lemma 2.1., we have \\

$~~~~~~ \begin{array}{l}
              < A_1 \eta, \chi > ~=~ < A \eta, \chi > - \overline{< A >_\psi}\, < \eta, \chi > ~=~ \\ \\
              ~=~ < \eta, A \chi > - < A >_\psi \,\, < \eta, \chi > ~=~
                        < \eta, A \chi > - < \eta, < A >_\psi\chi > ~=~ \\ \\
                       ~=~ < \eta, A_1 \chi >
        \end{array} $ \\

and similarly with $B_1$. \\

Next we prove \\

(2.19)~~~ $ [ A_1, B_1 ] = [ A, B ] $ \\

which is obtained as follows. We have from (2.17) \\

$~~~~~~ A_1 B_1 = A B \,\,- < A >_\psi B \,\,- < B >_\psi A \,\,+ < A >_\psi \,\,  < B >_\psi I $ \\

thus \\

$~~~~~~ B_1 A_1 = B A \,\,- < B >_\psi A \,\,- < A >_\psi B \,\,+ < B >_\psi \,\,  < A >_\psi I $ \\

hence (2.19). \\

Further, for $\psi \in {\cal S}_{\cal F} ( \mathbb{R} )$ with $ < \psi, \psi > \,= 1$, we have \\

(2.20)~~~ $ < A_1 \psi, A_1 \psi > ~=~ ( \Delta_{\psi} A )^2 $ \\

Indeed, in view of (2.18), we obtain \\

$~~~~~~ \begin{array}{l}
                < A_1 \psi, A_1 \psi > ~=~ < \psi, ( A_1 )^2 \psi > ~=~
                                       < \psi, ( A \,\,- < A >_\psi I )^2 \, \psi > ~=~ \\ \\
                 ~=~ < \psi, ( A^2 - 2 < A >_\psi A + ( < A >_\psi )^2 I ) \, \psi > ~=~ \\ \\
                 ~=~ < \psi, A^2 \psi > -\, 2 < A >_\psi \,\, < \psi, A \psi > + \,
                                     ( < A >_\psi )^2 < \psi, \psi > ~=~ \\ \\
                 ~=~ < A^2 >_\psi -\, 2 ( < A >_\psi )^2 + ( < A >_\psi )^2 ~=~ < A^2 >_\psi - ( < A >_\psi )^2
        \end{array} $ \\

In view of the above, we have for $\psi \in {\cal S}_{\cal F} ( \mathbb{R} )$ with $ < \psi, \psi > \,= 1$, the relations \\

$~~~~~~ \begin{array}{l}
            < \psi, [ A, B ] \psi > ~=~ < \psi, [ A_1, B_1 ] \psi > ~=~ < \psi, A_1 B_1 \psi > -
                                    < \psi, B_1 A_1 \psi > ~=~ \\ \\
            ~=~ < A_1 \psi, B_1 \psi > - < B_1 \psi, A_1 \psi > ~=~ < A_1 \psi, B_1 \psi > -
                                    \overline{< A_1 \psi, B_1 \psi >} ~=~ \\ \\
            ~=~ 2 i \, Im < A_1 \psi, B_1 \psi >
         \end{array} $ \\

Consequently \\

$~~~~~~ | \,  < \psi, [ A, B ] \psi >  \, | ~=~ 2 | \, Im < A_1 \psi, B_1 \psi > \, | $ \\

However, in view of (1.8), (1.9), (1.14), we have for $z \in \mathbb{C}_{\cal F}$ the relations \\

$~~~~~~ |\, Re \, z \,|,~~ |\, Im \, z \,| ~\leq~ |\, z \,| $ \\ \\

{\bf 3. The Extended Wintner Theorem on Unbounded \\
        \hspace*{0.5cm} Operators} \\

In the context of the Heisenberg uncertainty relation, the non-commutativity of operators involved is crucial, since for commutative operators the
respective inequality is obviously trivially satisfied, thus there is never any uncertainty, see (2.15). \\

Within the usual mathematical context of Quantum Mechanics given by complex Hilbert spaces $H$ and linear operators on them, it turns out that certain
simple non-commutativity relations for linear operators, such as (3.1) below, necessarily imply their unboundedness. Thus the need to consider unbounded,
however, densely defined and closed operators on such Hilbert spaces, and the fundamental operators of position and momentum are well known to be among
them. \\

Needless to say, this fact is a rather inconvenient one, since it complicates considerably the mathematical apparatus involved in Quantum Mechanics. And
it was precisely the avoidance of such a complication which led von Neumann to his second mathematical model for Quantum Mechanics, namely, the one based
on starting with algebras of observables, and then defining the states. Such an approach is obviosuly a reversal of the way in von Neumann's first model
based on Hilbert spaces of states, where the observables are then defined as Hermitian operators. \\

The classical result regarding the inevitability of unbounded operators in von Neumann's first model of Quantum Mechanics is given in \\

{\bf Wintner's Theorem} \\

Let $H$ be a complex Hilbert space and $A, B$ two bounded linear operators on it. Then there cannot be any nonzero constant $c \in \mathbb{C}$, such that
the non-commutation relation holds \\

(3.1)~~~ $ [ A, B ] = A B - B A = c I $ \\

where $I$ is the identity operator on $H$.

\hfill $\Box$ \\

The special relevance of this result is in the fact that the position and momentum operators do satisfy a non-commutation relation of type (3.1), this
therefore being the reason they cannot be given by bounded operators. \\

Here, we give an extended version of Wintner's Theorem to the case of linear operators in  ${\cal L} ( {\cal S}_{\cal F} ( \mathbb{R} ) )$, see (2.2),
(2.10), which are defined based on scalars given by reduced power algebras $\mathbb{R}_{\cal F}$ or $\mathbb{C}_{\cal F}$ . \\

The following simple linear functional analytic notions, extended to the case of reduced power algebras, will be needed. \\

For any wave function $\psi = \sum_{1 \leq h \leq m} \gamma_h H_h \in {\cal S}_{\cal F} ( \mathbb{R} )$, see (2.1), (2.2), we define its {\it norm} by \\

(3.2)~~~ $ ||\, \psi \,|| = \sup_{1 \leq h \leq m} |\, \gamma_h \,| \in \mathbb{R}^+_{\cal F} $ \\

We note that, with values in $\mathbb{R}^+_{\cal F}$, and not merely in $\mathbb{R}^+ = [ 0, \infty )$, as is the usual case with wave functions $\psi
\in L^2 ( \mathbb{R} )$, this norm (3.2) is always well defined, regardless of the $\gamma_h$ being infinitesimal, finite or infinitely large elements
in $\mathbb{C}_{\cal F}$. \\

Now, a linear operator $A \in {\cal L} ( {\cal S}_{\cal F} ( \mathbb{R} ) )$ is called {\it bounded}, if and only if there exists $M \in
\mathbb{R}^+_{\cal F}$, such that \\

(3.3)~~~ $ ||\, A \psi \,|| ~\leq~ M ||\, \psi \,||,~~~ \psi \in {\cal S}_{\cal F} ( \mathbb{R} ) $ \\

and in this case we denote by \\

(3.4)~~~ $ {\cal M}_A \subseteq \mathbb{R}^+_{\cal F} $ \\

the set of all such $M$. \\

We will also need the following partial order relation on $\mathbb{R}_{\cal F}$. Given $u = ( u_n )_{n \in \mathbb{N}} + {\cal I}_{\cal F},~ v =
( v_n )_{n \in \mathbb{N}} + {\cal I}_{\cal F} \in \mathbb{R}_{\cal F}$, we define \\

(3.5)~~~ $ u ~\leq~ v ~~\Longleftrightarrow~~ v - u \in \mathbb{R}^+_{\cal F} $ \\

Clearly, with the partial order (3.5), ${\cal M}_A$ has the property \\

(3.6)~~~ $ M\,' \in \mathbb{R}^+_{\cal F},~ M\,' \geq M \in {\cal M}_A ~~\Longrightarrow~~ M\,' \in {\cal M}_A $ \\

The interest in dealing with ${\cal M}_A$ is that, in this way, we can avoid the issue of considering the existence, and
of the properties of $\inf ~{\cal M}_A$ in the reduced power algebra $\mathbb{R}_{\cal F}$. Here we note that in the
usual case of operators on Hilbert spaces, instead of (3.4) one obviously has ${\cal M}_A \subseteq \mathbb{R}^+ = [ 0,
\infty )$, thus  $\inf ~{\cal M}_A \in \mathbb{R}^+$ always exists, and it is denoted by $|| A ||$, hence the above issue
simply does not arise. \\
However, as the following two easy to prove Lemmas show it, we can to a good extent avoid that issue even in the general
case of arbitrary reduced power algebras $\mathbb{R}_{\cal F}$. \\

{\bf Lemma 3.1.} \\

Let $\psi,~ \psi\,' \in {\cal S}_{\cal F} ( \mathbb{R} )$ and $c \in {\mathbb{C}_{\cal F}}$, then \\

~~~~~~ $ 1)~~ ||\, \psi + \psi\,' \,|| \leq  ||\, \psi \,|| + ||\, \psi\,' \,|| $ \\

~~~~~~ $ 2)~~ ||\, c \, \psi \,|| = |\, c \,| \,\, ||\, \psi \,|| $ \\

~~~~~~ $ 3)~~ ||\, \psi \,|| = 0 ~~~ \Longleftrightarrow~~~ \psi = 0 \in {\cal S}_{\cal F} ( \mathbb{R} ) $ \\

{\bf Lemma 3.2.} \\

 Let be any bounded operators $A, B \in {\cal L} ( {\cal S}_{\cal F} ( \mathbb{R} ) )$ and $c \in \mathbb{C}_{\cal F}$. Then the following hold \\

 $ \begin{array}{l}
             1)~~~  \forall~~ K \in  {\cal M}_A,~ L \in {\cal M}_B ~: \\ \\
                 ~~~~~~ \exists~~ M \in {\cal M}_{A + B} ~: \\ \\
                 ~~~~~~~~~~ M \leq K + L
      \end{array} $ \\ \\

 $ \begin{array}{l}
             2)~~~  \forall~~ K \in  {\cal M}_A ~: \\ \\
                 ~~~~~~ \exists~~ M \in {\cal M}_{c A} ~: \\ \\
                 ~~~~~~~~~~ M \leq | c | K
      \end{array} $ \\ \\

$ \begin{array}{l}
             3)~~~  \forall~~ K \in  {\cal M}_A,~ L \in {\cal M}_B ~: \\ \\
                 ~~~~~~ \exists~~ M \in {\cal M}_{A B} ~: \\ \\
                 ~~~~~~~~~~ M \leq K L
      \end{array} $ \\ \\

The extension of the Wintner Theorem can now be formulated as follows \\

{\bf Theorem 3.1.} \\

Let $A, B \in {\cal L} ( {\cal S}_{\cal F} ( \mathbb{R} ) )$ two bounded linear operators. Then there cannot be any nonzero constant $c \in
\mathbb{C}_{\cal F}$, such that the non-commutation relation holds \\

(3.7)~~~ $ [ A, B ] = A B - B A = c I $ \\

where $I$ is the identity operator on ${\cal S}_{\cal F} ( \mathbb{R} )$. \\

{\bf Proof.} \\

Obviously, it suffices to consider the case $c = 1$. Let us then assume $A, B \in {\cal L} ( {\cal S}_{\cal F}
( \mathbb{R} ) )$ two bounded linear operators, such that \\

(3.8)~~~ $ A B - B A = I $ \\

Then by induction, one obtains \\

(3.9)~~~ $ n B^{n - 1} = A B^n - B^n A,~~~ n \geq 1 $ \\

Indeed, for $n = 1$, the relation (3.9) reduces to (3.8). Assuming now that (3.9) holds for a certain $n \geq 1$, we have
then \\

(3.10)~~~ $ \begin{array}{l}
                  ( n + 1 ) B^n ~=~ n B^{n - 1} B + B^n I ~=~ \\ \\
                  ~~~=~ ( A B^n - B^n A ) B + B^n ( A B - B A ) + A B^{n + 1} - B^{n + 1} A
             \end{array} $ \\

Now in view of Lemma 3.2. applied to (3.9), one obtains \\

(3.11)~~~ $ \begin{array}{l}
                      \forall~~ n \geq 1,~ K \in  {\cal M}_A,~ L \in {\cal M}_B ~: \\ \\
                      \exists~~ M \in {\cal M}_{B^{n - 1}} ~: \\ \\
                      ~~~~ n M \leq 2 K L M
             \end{array} $ \\

consequently, we must have $M = 0$, for some $n \geq 1$. And then, (3.9) implies $B^{n - 1} = 0$, and thus successively,
$B = 0$, and finally $I = 0$, which of course is absurd. \\ \\

{\bf 4. Question on Two Fundamental Physical Constants} \\

Re-formulating Physics in terms of scalars given by reduced power algebras leads naturally to the following two questions

\begin{itemize}

\item Is it possible that Planck's constant $h$ is in fact an {\it infinitesimal} in some reduced power algebra $\mathbb{R}_{\cal F}$ ?

\item Is it possible that the maximum speed of propagation of physical effects is not finite, but rather an {\it infinitely large} quantity in some
reduced power algebra $\mathbb{R}_{\cal F}$ ?

\end{itemize}

The motivation for these two questions appears quite natural, as
soon as one becomes more familiar with the non-Archimedean structure
of reduced power algebras, [2-9]. Indeed, that non-Archimedean
structure leads to the presence of three types of elements in such
algebras, namely : infinitesimals,
finite elements, and infinitely large elements. \\

The essential fact in this regard, however, is that the above classification in three types of elements is {\it relative}. Namely, it is implied by the
fact that, when constructing reduced power algebras $\mathbb{R}_{\cal F}$, one starts by defining the usual real numbers in $\mathbb{R}$ as being the
finite ones. Indeed, such reduced power algebras have a highly complex and rich {\it self-similar} structure. And it is easy to see that, due to that
structure, one is in fact {\it not} obliged to choose the usual real numbers in $\mathbb{R}$ as being the finite ones. On the contrary, that
self-similar structure renders the concept of "finite elements" highly relative, by allowing a wide range of other choices. In this way, elements which
in a choice are finite, may become infinitesimal or infinitely large in other choices, and vice-versa. \\

Consequently, when using reduced power algebras in Physics, thus non-Archimedean scalar structures, one is {\it no longer} obliged to have both the
Planck constant and that of the maximum speed of propagation of physical effects finite, and thus having the only possible difference between between
them reduced to a large but finite factor. \\

\end{document}